Article

# Stochastic Multi-Objective Multi-Trip AMR Routing Problem with Time Windows

Lulu Cheng [1], Ning Zhao [1] and Kan Wu [2,*]

[1] Faculty of Science, Kunming University of Science and Technology, Kunming 650500, China
[2] Business Analytics Research Center, Chang Gung University, Taoyuan City 33302, Taiwan
* Correspondence: kan626@gmail.com

**Abstract:** In recent years, with the rapidly aging population, alleviating the pressure on medical staff has become a critical issue. To improve the work efficiency of medical staff and reduce the risk of infection, we consider the multi-trip autonomous mobile robot (AMR) routing problem in a stochastic environment. Our goal is to minimize the total expected operating cost and maximize the total service quality for patients, ensuring that each route violates the vehicle capacity and the time window with only a minimal probability. The travel time of AMRs is stochastically affected by the surrounding environment; the demand for each ward is unknown until the AMR reaches the ward, and the service time is linearly related to the actual demand. We developed a population-based tabu search algorithm (PTS) that combines the genetic algorithm with the tabu search algorithm to solve this problem. Extensive numerical experiments were conducted on the modified Solomon instances to demonstrate the efficiency of the PTS algorithm and reveal the impacts of the confidence level on the optimal solution, providing insights for decision-makers to devise delivery schemes that balance operating costs with patient satisfaction.

## 1. Introduction

In recent years, with the increasing global aging population, medical services have become a major challenge [1]. This aging trend has led to a rise in chronic diseases and a growing demand for healthcare resources. Consequently, there is a shortage of medical staff, which further exacerbates the pressure and burden on the medical industry [2]. Transportation tasks in nursing directly impact patient care time in the context of staff shortages and an aging society [3]. In addition, healthcare workers face a higher risk of infection when dealing with infectious diseases or outbreaks, which not only threatens their health but may also affect the operation and service quality of medical institutions.

With the rapid development of the smart manufacturing industry, robots can take over routine transportation tasks, such as lab specimens and patient forms, which ties up nursing staff and detracts from time devoted to patient care, treatment, and therapy [4]. Currently, autonomous mobile robots (AMRs) have been widely used in various fields, such as material handling, unmanned aerial vehicle cruises, automatic factories, and so forth [5]. An increasing number of scholars are dedicating their efforts to investigating AMR and automated guided vehicle (AGV) scheduling problems [6–8]. In real application scenarios, stochastic factors play a vital role in the routing and scheduling of AMRs. Sources of stochastic in an AMR scheduling problem can be travel times, service times, and demands of customers. Ignoring these stochastic factors may result in AMR scheduling problems that are not applicable in real situations. Aziez et al. [9] studied the fleet sizing and routing problem with synchronization for AGVs with dynamic demands in the context of a real-life application and proposed a powerful metaheuristic based on a fast and efficient dynamic re-optimization of the routes upon the arrival of new requests. Cheng et al. [10] investigated the optimization problem of scheduling strategies for AMRs at smart hospitals, where the service and travel times of AMRs are stochastic, but the demand for each medical request is deterministic.

However, the factors considered in these studies are not comprehensive enough. Primarily, AMR scheduling is subject to various uncertain factors during operation, resulting in randomness. Furthermore, there is uncertainty in the demand for each request, and the service time depends on different demand levels. Lastly, the AMR exhibits a limited capacity, and each AMR will travel multiple trips to complete the service. Therefore, to address the above issues practically, we broadened previous research to study scheduling problems that are more in line with practical applications. This paper investigates a stochastic multi-objective, multi-trip autonomous mobile robot routing problem with time windows (SMMAMRRPTW). This problem involves multiple homogeneous AMRs assisting medical staff in completing related delivery work. The travel time of AMRs is stochastic, and the demand for each ward is uncertain. Additionally, the service time of AMRs depends on the size of the demand for the ward. Each AMR travels multiple paths to provide services to various wards while meeting its payload. The objective of the problem is to minimize the cost of investing AMRs in hospitals and maximize patient service efficiency. The main contributions of this paper are as follows:

1. Formulate the SMMAMRRPTW and establish a stochastic programming model.;
2. Develop a population-based tabu search algorithm (PTS);
3. Demonstrate the effectiveness and robustness of the PTS algorithm through numerical experiments, compared with other algorithms, and provide relevant management insights.

The rest of the paper is as follows: Section 2 reviews the literature on AMR scheduling problems and the stochastic vehicle routing problem. Section 3 introduces the stochastic programming model for solving the problem. Section 4 describes the proposed population-based tabu search algorithm. Section 5 verifies the effectiveness of the proposed algorithm by a large number of experiments using a modified Solomon instance. Finally, the conclusions and looks forward to the future follow in Section 6.

## 2. Literature Review

Recently, an increasing number of researchers have dedicated their efforts to investigating AMR scheduling problems due to the development of the smart manufacturing industry. However, there is relatively little literature on the AMR scheduling problem in a stochastic environment. Therefore, the following literature review focuses on the AMR or AGV scheduling problem and stochastic vehicle routing problems.

As for the closely related AGV or AMR scheduling problem in flexible manufacturing systems and the healthcare industry, Nishi et al. [11] addressed a bilevel decomposition algorithm to solve the simultaneous scheduling and conflict-free routing problems for AGVs. Chikul et al. [12] developed a deterministic model for hospital supply chain management, evaluating the differences between manual and technologically integrated approaches. Zou et al. [13] addressed a new multiple AGV dispatching problem. A mixed integer linear programming model was formulated to minimize transportation costs, and a discrete artificial bee colony algorithm was presented to solve the problem. Han et al. [14] proposed a digital twin-based dynamic AGV scheduling method to solve the AGV charging problem in the AGV scheduling system. Cheng et al. [15] aimed to minimize daily costs in hospitals and proposed an improved tabu search algorithm to solve the AMR scheduling problem with known medical request demand. However, due to the distinct nature of our specific research objectives, direct application of the approaches above may prove challenging.

SMMAMRRPTW shares many characteristics with stochastic vehicle routing problems and is a variant of stochastic vehicle routing problems. Gendreau et al. [16] proposed the real-time dynamic vehicle dispatching problem with pick-up and deliveries for new requests (including a pick-up and a delivery location). Biesinger et al. [17] considered the generalized vehicle routing problem with stochastic demands. Goel et al. [18] considered vehicle routing problems with time windows having stochastic demands and stochastic service times. Çimen and Soysal [19] addressed a time-dependent capacitated vehicle routing problem with stochastic vehicle speeds and environmental concerns. This problem takes time-dependency and stochasticity of the vehicle speeds into account. Elgharably et al. [20] studied the stochastic multi-objective vehicle routing problem in a green environment, where stochastic factors include demand type, service time, and travel time. Xie et al. [21] studied the dynamic takeout delivery vehicle routing problem, which faces fluctuating order demand and time-varying speeds, and proposed a dynamic takeout delivery vehicle routing optimization model. The research on the vehicle routing problem mentioned above mostly considers incomplete random factors and neglects scenarios involving multi-trip transportation. Thus, addressing the SMMAMRRPTW is crucial for practical applications.

The vehicle routing problem has proven to be NP-hard in the strong sense [22]. In general, it is almost impossible for an NP-hard problem to obtain a solution by using exact solution methods in a limited amount of computing time [23]. For an NP-hard problem, heuristic and metaheuristic approaches are more suitable for solving more complex and large instances than the exact methods. Gendreau et al. [16] proposed a neighborhood search heuristic algorithm to optimize the planned route of vehicles in real-time when new requests occur. Biesinger et al. [17] presented a steady-state genetic algorithm to solve the generalized vehicle routing problem with stochastic demands using a preventive restocking strategy. Goel et al. [18] proposed a modified ant colony system to solve the vehicle routing problem with time windows having stochastic demands and stochastic service times to minimize total transportation costs and maximize customer satisfaction. Elgharably et al. [20] presented a new hybrid search algorithm to solve the stochastic multi-objective vehicle routing problem. Li and Li [24] proposed an improved tabu search algorithm based on a greedy algorithm to solve the vehicle routing problem with stochastic travel time and service time. Wu et al. [25] designed a tailored iterative local search heuristic for solving the multi-trip vehicle routing problem with multiple time windows. And the effectiveness of the proposed model and algorithm was experimentally demonstrated. The aforementioned literature has algorithms for the study of vehicle routing problems in a stochastic environment; however, the proposed algorithm cannot be directly applied to the SMMAMRRPTW.

As previously mentioned, the scheduling problem in a stochastic environment has become a research hotspot in recent years. However, no research on SMMAMRRPTW has been documented in the existing literature. Even though existing mathematical programming model and optimization algorithms have been proposed for scheduling problems in a stochastic environment, they cannot be used directly for the SMMAMRRPTW due to their problem-specific characteristics. Therefore, SMMAMRRPTW must establish a mathematical programming model and propose an efficient algorithm is vital for SMMAMRRPTW.

## 3. Problem Description and Formulation

The characteristics of the proposed stochastic multi-objective multi-trip autonomous mobile robot routing problem with time windows (SMMAMRRPTW) are presented, followed by the mathematical formulation of the problem.

*3.1. Problem Description*

In real life, each ward generates an unknown amount of medical waste, which is transported by AMRs to assist staff, and the actual demand only becomes revealed when AMRs reach the ward. Assuming the hospital introduces a fleet of $m$ homogeneous AMRs with capacity $Q$ to participate in medical services, which is stationed at the depot (waste disposal point), they will be waiting for the medical staff's instructions to pick up the medical waste to $n$ medical requests (wards). It is worth noting that the number of AMRs considered in this paper is not fixed but is to be determined by the solution approach. To improve the service rate of AMRs and reduce the cost of investing in AMRs, this research allows each AMR to travel multiple trips, and the set of routes traveled by the $k$th AMR is denoted as $L_k = \{1, 2, \cdots, l_{max}^k\}$. AMRs can avoid obstacles and take the elevator when serving medical requests. The uncertainty of the road conditions leads to the random nature of the travel time.

The following is a complete graph $G = (V, A)$ to define the SMMAMRRPTW, where $V = \{0, 0'\} \cup R$ is the vertex set, and $A = \{(i,j) | i, j \in V, i \neq j\}$ is the edge set. The vertex $0'$ is a dummy copy of the depot $0$ with the same location. $R = \{1, 2, \cdots, n\}$ is the set of medical requests. Each request $i \in R$ has its time window $[e_i, h_i]$, service time $S_i$, and demand $q_i$, for each AMR, its load cannot exceed a given boundary $Q$. When the $k$th AMR arrives at request $i \in R$ on the $p \in L_k$ path, the time is $A_{pi}^k$, the waiting time is $W_{pi}^k$, the time to start service request $i$ is $Y_{pi}^k$, the remaining capacity of the AMR is $u_{pi}^k$. The travel time of the $k$th AMR from request $i \in R$ to request $j \in R$ on its $p$th route is $T_{pij}^k$. AMRs are required to provide services for requests within a prespecified time window. The uncertainty in travel time and service time would make the on-time service a challenging one to enforce. Hence, if an AMR arrives before or after the time window, it has a penalty cost.

It was found that most studies are based on the Central Limit Theorem, which states that the distribution of sufficiently large random samples will be approximately normally distributed. Assuming that the uncertain parameters follow the normal distribution [26]. Correspondingly, we assume that the demands ($q_i, i \in R$) are independent and follow a normal distribution with mean $\mu_{q_i}$ and variance $\sigma_{q_i}^2$, the travel time follows a normal distribution with mean $\mu_{T_{pij}^k}$ and variance $\sigma_{T_{pij}^k}^2$. Since the actual service time of each ward is affected by the actual demand, we use $S_i = a \cdot q_i + b$ to represent the service time of request $i$, where $a$ is the service time per demand unit and $b$ is a bias. Hence, the service time also follows a normal distribution with mean $\mu_{S_i} = a \cdot \mu_{q_i} + b$ and variance $\sigma_{S_i}^2 = a^2 \cdot \sigma_{q_i}^2$.

Table 1 provides a summary of the notation used. The parameters and decision variables are defined as follows.

**Table 1.** Notation for the SMMAMRRPTW.

| | |
|---|---|
| **Sets** | |
| $R = \{1, 2, ..., n\}$ | Set of requests |
| $L_k$ | Set of the route traveled by the $k$th AMR |
| **Parameters** | |
| $Q$ | Cargo capacity of the AMR |
| $\xi_1$ | Fixed cost of the AMR |
| $\xi_2$ | Cost per unit time |
| $q_i$ | Demand of request $i$ |

| $S_i$ | Service time of request $i$ |
|---|---|
| $W_{pi}^k$ | Waiting time of the $k$th AMR at the request $i$ on the $p$th route |
| $T_{pij}^k$ | Travel time spent by the $k$th AMR on the route $p \in L_k$ from request $i$ to request $j$ |
| $u_{pi}^k$ | Remaining capacity of the $k$th AMR when it arrives at request $i$ on its $p$th route |
| $A_{pi}^k$ | Time when the $k$th AMR arrives at the request $i$ on the $p$th route |
| $Y_{pi}^k$ | Start service request $i$ time of the $k$th AMR on the $p$th route |
| $[e_i, h_i]$ | Time window of request $i$ |
| **Decision variables** | |
| $x_{pij}^k$ | 1 if the $k$th AMR completes request $i$ and $j$ successively on its $p$th route; 0 otherwise |
| $m$ | Number of AMRs participating in medical request services |
| $z_{pi}^k$ | 1 if the $k$th AMR completes service request $i \in R$ on its $p$th route; 0 otherwise |

### 3.2. Model Formulation

With the above description, we present a mathematical model for the SMMAMRRPTW. Chance-constrained programming (CCP) is a modeling method for vehicle routing problems with stochastic demands. In the CCP, route failure is accepted with a probability of failure, no corrective action is taken to satisfy customer demands, and the costs of those failures are ignored [27]. This section establishes a chance-constrained stochastic multi-objective AMR routing programming model to deal with demand and time uncertainties. We aim to minimize the sum of fixed and time costs of AMR and maximize the quality of patient services.

The problem is formulated as follows:

$$\text{minmize } f_1 = \xi_1 m + \xi_2 \sum_{p \in L_k} \sum_{k=1}^{m} \sum_{j \in V} \sum_{i \in V} E\left[\left(T_{pij}^k + S_i + W_{pi}^k\right) \cdot x_{pij}^k\right] \quad (1)$$

$$\text{maxmize } f_2 = \sum_{p \in L_k} \sum_{k=1}^{m} \sum_{i \in R} E\left[\left(h_i - A_{pi}^k\right) \cdot z_{pi}^k\right]^+ \quad (2)$$

subject to:

$$\sum_{j \in R \cup \{0'\}, j \neq i} x_{pij}^k = \sum_{j \in \{0\} \cup R, j \neq i} x_{pji}^k, i \in V, k = 1, \ldots, m, p \in L_k \quad (3)$$

$$\sum_{p \in L_k} \sum_{k=1}^{m} \sum_{j \in V \setminus \{0\}, i \neq j} x_{pij}^k = 1, i \in R \quad (4)$$

$$\sum_{i \in V \setminus \{0'\}} x_{pij}^k = z_{pj}^k, j \in V \setminus \{0\}, k = 1, \ldots, m, p \in L_k \quad (5)$$

$$u_{pj}^k \leq u_{pi}^k - x_{pij}^k \cdot q_{pi}^k + Q(1 - x_{pij}^k), i \in V \setminus \{0'\}, j \in V \setminus \{0\}, i \neq j, k = 1, \ldots, m, p \in L_k \quad (6)$$

$$P\left\{\sum_{i \in R} q_i \cdot z_{pi}^k \leq Q\right\} \geq 1 - \varepsilon_1, k = 1, \ldots, m, p \in L_k \quad (7)$$

$$P\left\{q_j \cdot x_{pij}^k \leq u_{pi}^k\right\} \geq 1 - \varepsilon_1, i \in V \setminus \{0'\}, j \in V \setminus \{0\}, i \neq j, k = 1, \ldots, m, p \in L_k \quad (8)$$

$$P\left\{Y_{pi}^k + S_i + T_{pij}^k \leq h_j\right\} \geq 1 - \varepsilon_2, i, j \in R, i \neq j, k = 1, \ldots, m, p \in L_k \quad (9)$$

$$Y_{pi}^k + S_i + T_{pij}^k \cdot x_{pij}^k - H\left(1 - x_{pij}^k\right) \leq Y_{pj}^k, i \in V \setminus \{0'\}, j \in V \setminus \{0\}, i \neq j, k = 1, \ldots, m, p \in L_k \quad (10)$$

$$Y_{pi}^k = A_{pi}^k + W_{pi}^k, i \in R, k = 1,\ldots,m, p \in L_k \tag{11}$$

$$Y_{pi}^k = Y_{p+1,j}^k, i \in \{0'\}, j \in \{0\}, k = 1,\ldots,m, p \in L_k \tag{12}$$

$$W_{pi}^k = \max\{0, e_i - A_{pi}^k\}, i \in V, k = 1,\ldots,m, p \in L_k \tag{13}$$

$$x_{pij}^k \in \{0,1\}, i \in V \setminus \{0'\}, j \in V \setminus \{0\}, i \neq j, k = 1,\ldots,m, p \in L_k \tag{14}$$

$$z_{pi}^k \in \{0,1\}, i \in R, k = 1,\ldots,m, p \in L_k \tag{15}$$

Equation (1) represents the first objective function of minimizing the sum of fixed cost and time cost of AMRs. Equation (2) represents the second objective function of maximizing the quality of service to patients. Here, we assume that all AMRs have an identical fixed cost $\xi_1$ (\$/each) and unit time cost $\xi_2$ (\$/sec); meanwhile, $\xi_1 > \xi_2$. Equation (3) indicates that the number of input arcs at each point is equal to the number of output arcs. Equations (4) and (5) impose that each request must be serviced exactly once by one route. Inequalities (6), (7), and (8) represent the capacity constraints of the AMR on any route. Inequalities (7) and (8) ensure that the probability of satisfying the capacity of the AMRs is not less than the confidence level $1-\varepsilon_1$. The inequality (9) ensures that request $i$ is served in its time window with a confidence level of $1-\varepsilon_2$. The constraint (10) represents the relationship of the start service time between a task and its predecessor, where $H$ is a large number. Equations (11)–(13) represent the expressions for the start service time and waiting time of AMR. Formulas (14) and (15) denote binary decision variables $x_{pij}^k$ and $z_{pi}^k$.

Based on the start service time, travel time, and service time, we can express the time when the $p$th path of the $k$th AMR reaches request $i$ as $A_{pi}^k = Y_{pi-1}^k + S_{i-1} + T_{p,i-1,i}^k$. We take the approach of Ehmke et al. [28] to approximate the arrival time by a normal distribution, i.e., $A_{pi}^k \sim N\left(\mu_{A_{pi}^k}, \sigma_{A_{pi}^k}^2\right)$, where $\mu_{A_{pi}^k} = \mu_{Y_{pi-1}^k} + \mu_{S_{i-1}} + \mu_{T_{pi-1,i}^k}$ and $\sigma_{A_{pi}^k}^2 = \sigma_{Y_{pi-1}^k}^2 + \sigma_{S_{i-1}}^2 + \sigma_{T_{pi-1,i}^k}^2$. Therefore, the objective function (2) can be expressed as

$$\sum_{p \in L_k}\sum_{k=1}^m\sum_{i \in R} E\left[(h_i - A_{pi}^k) \cdot z_{pi}^k\right]^+ = \sum_{p \in L_k}\sum_{k=1}^m\sum_{i \in R} E\left(\max\{h_i, A_{pi}^k\} - A_{pi}^k\right) \cdot z_{pi}^k$$

$$= \sum_{p \in L_k}\sum_{k=1}^m\sum_{i \in R}\left(\begin{array}{l} h_i \Phi\left(\dfrac{h_i - \mu_{A_{pi}^k}}{\sigma_{A_{pi}^k}}\right) + \mu(A_{pi}^k)\Phi\left(\dfrac{\mu_{A_{pi}^k} - h_i}{\sigma_{A_{pi}^k}}\right) \\ + \sigma_{A_{pi}^k}\phi\left(\dfrac{h_i - \mu_{A_{pi}^k}}{\sigma_{A_{pi}^k}}\right) - \mu_{A_{pi}^k} \end{array}\right) \cdot z_{pi}^k. \tag{16}$$

The waiting time $W_{pi}^k = \max\{0, e_i - A_{pi}^k\}$ is a random variable because it depends on travel time $T_{pij}^k$, service time $S_i$, and the waiting time of the previous node. Nadarajah and Kotz [29] studied the mean and variance of the maximum of two normally distributed variables, $X = \max\{X_1, X_2\}$ with means $\mu_1$ and $\mu_2$ and variances $\sigma_1^2$ and $\sigma_2^2$, and correlation coefficient $\rho$ are known. Therefore, the expected and variance of the $W_{pi}^k$ are

$$\mu_{W_{pi}^k} = \mu_{A_{pi}^k}\Phi\left(\dfrac{\mu_{A_{pi}^k} - e_i}{\sigma_{A_{pi}^k}}\right) + e_i\Phi\left(\dfrac{e_i - \mu_{A_{pi}^k}}{\sigma_{A_{pi}^k}}\right) + \sigma_{A_{pi}^k}\phi\left(\dfrac{\mu_{A_{pi}^k} - e_i}{\sigma_{A_{pi}^k}}\right) - \mu_{A_{pi}^k} \tag{17}$$

$$\sigma^2_{W^k_{pi}} = \left(\mu^2_{A^k_{pi}} + \sigma^2_{A^k_{pi}}\right)\Phi\left(\frac{\mu_{A^k_{pi}} - e_i}{\sigma_{A^k_{pi}}}\right) + e_i^2 \Phi\left(\frac{e_i - \mu_{A^k_{pi}}}{\sigma_{A^k_{pi}}}\right) - \mu^2_{W^k_{pi}} - \sigma^2_{A^k_{pi}}$$
$$+ \left(e_i + \mu_{A^k_{pi}}\right)\sigma_{A^k_{pi}} \phi\left(\frac{\mu_{A^k_{pi}} - e_i}{\sigma_{A^k_{pi}}}\right) \quad (18)$$

$\Phi(\cdot)$ is the cumulative distribution function of the standard normal distribution, and $\phi(\cdot)$ is the probability density function of the standard normal distribution. Therefore, the objective function (1) can be expressed as

$$\xi_1 m + \xi_2 \sum_{p \in L_k} \sum_{k=1}^{m} \sum_{j \in V} \sum_{i \in V} \left[\left(\mu_{T^k_{pij}} + \mu_{S_i} + \mu_{W^k_{pi}}\right) \cdot x^k_{pij}\right]. \quad (19)$$

Based on the above discussion, the objective of the model can be transformed into a deterministic expression composed of Equations (16) and (19). And according to Equation (13), we also know that the expected start time of the service is $\mu_{Y^k_{pi}} = \mu_{A^k_{pi}} + \mu_{W^k_{pi}}$, and the variance is $\sigma^2_{Y^k_{pi}} = \sigma^2_{A^k_{pi}} + \sigma^2_{W^k_{pi}}$.

Remaining capacity $u^k_{pi}$ is a random variable because it depends on demand $q_i$ which is defined as a random variable. The following lemmas explain the distribution and calculation of (7)–(9). Below, the reformulation of inequalities (7)–(9) is explained in the following to find a deterministic optimization model.

**Lemma 1.** *If $\tilde{q}^k_p = \sum_{i \in R} q_i z^k_{pi}$ and $z^k_{pi}$ is a binary variable, then we have:*

$$P\left\{\sum_{i \in R} q_i \cdot z^k_{pi} \leq Q\right\} \geq 1 - \varepsilon_1 \Leftrightarrow \Phi\left(\frac{Q - \sum_{i \in R} z^k_{pi} \mu_{q_i}}{\sqrt{\sum_{i \in R} z^k_{pi} \sigma^2_{q_i}}}\right) \geq 1 - \varepsilon_1, \forall k = 1, \dots, m,\ p \in L_k \quad (20)$$

**Proof.** Let $\tilde{q}^k_p = \sum_{i \in R} q_i z^k_{pi}$. Due to $\tilde{q}^k_p = \sum_{i \in R} q_i z^k_{pi}$ and $q_i \sim N\left(\mu_{q_i}, \sigma^2_{q_i}\right)$, then $\tilde{q}^k_p \sim N\left(\mu_{\tilde{q}^k_p}, \sigma^2_{\tilde{q}^k_p}\right)$, where

$$\mu_{\tilde{q}^k_p} = \sum_{i \in R} \mu_{q_i} z^k_{pi}, \forall k = 1, \dots, m,\ p \in L_k \quad (21)$$

$$\sigma^2_{\tilde{q}^k_p} = \sum_{i \in R} \sigma^2_{q_i} z^k_{pi}, \forall k = 1, \dots, m,\ p \in L_k \quad (22)$$

The normal random variable $\tilde{q}^k_p$ is converted to the standard normal distribution by transformation function $X_i = \dfrac{\tilde{q}^k_p - \mu_{\tilde{q}^k_p}}{\sigma_{\tilde{q}^k_p}}$. Let $p_1 = P\left\{\sum_{i \in R} q_i \cdot z^k_{pi} \leq Q\right\}$. So, the probability $p_1$ is

$$p_1 = P\left\{\tilde{q}^k_p \leq Q\right\} = P\left\{X_i \leq \frac{Q - \mu_{\tilde{q}^k_p}}{\sigma_{\tilde{q}^k_p}}\right\} = \Phi\left(\frac{Q - \mu_{\tilde{q}^k_p}}{\sigma_{\tilde{q}^k_p}}\right).$$

From Equations (21) and (22), we have

$$p_1 = \Phi\left(\frac{Q - \sum_{i \in R} z_{pi}^k \mu_{q_i}}{\sqrt{\sum_{i \in R} z_{pi}^k \sigma_{q_i}^2}}\right),$$

Hence,

$$P\left\{\sum_{i \in R} q_i \cdot z_{pi}^k \leq Q\right\} \geq 1 - \varepsilon_1 \Leftrightarrow \Phi\left(\frac{Q - \sum_{i \in R} z_{pi}^k \mu_{q_i}}{\sqrt{\sum_{i \in R} z_{pi}^k \sigma_{q_i}^2}}\right) \geq 1 - \varepsilon_1, \forall k = 1, \ldots, m, \ p \in L_k.$$

□

**Lemma 2.** *If* $q_i \sim N(\mu_{q_i}, \sigma_{q_i}^2)$, $z_{pi}^k$ *is a binary variable, and* $Q$ *is a constant, then we have*

$$P\{q_j \cdot x_{pij}^k \leq u_{pi}^k\} \geq 1 - \varepsilon_1 \Leftrightarrow \Phi\left(\frac{Q - \mu_{q_j} x_{pij}^k - \sum_{i \in R} \mu_{q_i} z_{pi}^k}{\sqrt{\sigma_{q_j}^2 x_{pij}^k + \sum_{i \in R} \sigma_{q_i}^2 z_{pi}^k}}\right) \geq 1 - \varepsilon_1, \tag{23}$$

$$\forall i \in \{0\} \cup R, j \in R \cup \{0'\}, i \neq j, k = 1, \ldots, m, p \in L_k$$

**Proof.** Because $q_i \sim N(\mu_{q_i}, \sigma_{q_i}^2)$ and $u_{pi}^k = Q - \sum_{i \in R} q_i z_{pi}^k$, we have $u_{pi}^k \sim N(\mu_{u_{pi}^k}, \sigma_{u_{pi}^k}^2)$, where

$$\mu_{u_{pi}^k} = Q - \sum_{i \in R} \mu_{q_i} z_{pi}^k, \forall k = 1, \ldots, m, \ p \in L_k \tag{24}$$

$$\sigma_{u_{pi}^k}^2 = \sum_{i \in R} \sigma_{q_i}^2 z_{pi}^k, \forall k = 1, \ldots, m, \ p \in L_k \tag{25}$$

Let $X_i = q_j x_{pij}^k - u_{pi}^k$; obviously, $X_i$ is a normal random variable with mean $\mu_{X_i} = \mu_{q_j} x_{pij}^k - \mu_{u_{pi}^k}$ and variance $\sigma_{X_i}^2 = \sigma_{q_j}^2 x_{pij}^k + \sigma_{u_{pi}^k}^2$. The normal random variable $X_i$ is converted to the standard normal distribution by transformation function $X_i' = \frac{X_i - \mu_{X_i}}{\sigma_{X_i}}$. Let $p_2 = P\{q_j x_{pij}^k \leq u_{pi}^k\} = P\{q_j x_{pij}^k - u_{pi}^k \leq 0\} = P\{X_i \leq 0\}$, So, the probability $p_2$ is

$$p_2 = P\{X_i \leq 0\} = P\left\{X_i' \leq \frac{-\mu_{X_i}}{\sigma_{X_i}}\right\} = \Phi\left(-\frac{\mu_{X_i}}{\sigma_{X_i}}\right) = \Phi\left(\frac{\mu_{u_{pi}^k} - \mu_{q_j} x_{pij}^k}{\sqrt{\sigma_{q_j}^2 x_{pij}^k + \sigma_{u_{pi}^k}^2}}\right).$$

From Equations (24) and (25), we have

$$p_2 = \Phi\left(\frac{Q - \mu_{q_j} x_{pij}^k - \sum_{i \in R} \mu_{q_i} z_{pi}^k}{\sqrt{\sigma_{q_j}^2 x_{pij}^k + \sum_{i \in R} \sigma_{q_i}^2 z_{pi}^k}}\right).$$

Hence,

$$P\{q_j \cdot x_{pij}^k \leq u_{pi}^k\} \geq 1 - \varepsilon_1 \Leftrightarrow \Phi\left(\frac{Q - \mu_{q_j} x_{pij}^k - \sum_{i \in R} \mu_{q_i} z_{pi}^k}{\sqrt{\sigma_{q_j}^2 x_{pij}^k + \sum_{i \in R} \sigma_{q_i}^2 z_{pi}^k}}\right) \geq 1 - \varepsilon_1,.$$

$$\forall i \in \{0\} \cup R, j \in R \cup \{0'\}, i \neq j, k = 1, \ldots, m, p \in L_k$$

□

**Lemma 3.** *If* $T_{pij}^k \sim N\left(\mu_{T_{pij}^k}, \sigma_{T_{pij}^k}^2\right)$, $A_{pi}^k \sim N\left(\mu_{A_{pi}^k}, \sigma_{A_{pi}^k}^2\right)$, $S_i \sim N\left(\mu_{S_i}, \sigma_{S_i}^2\right)$ *and* $W_{pi}^k$ *with mean Equation (17) and variance Equation (18), then we have:*

$$P\left\{Y_{pi}^k + S_i + T_{pij}^k \leq h_j\right\} \geq 1 - \varepsilon_2 \Leftrightarrow \Phi\left(\frac{h_j - \mu_{A_{pi}^k} - \mu_{W_{pi}^k} - \mu_{S_i} - \mu_{T_{pij}^k}}{\sqrt{\sigma_{A_{pi}^k}^2 + \sigma_{W_{pi}^k}^2 + \sigma_{S_i}^2 + \sigma_{T_{pij}^k}^2}}\right) \geq 1 - \varepsilon_2, \quad (26)$$

$$i, j \in R, i \neq j, k = 1, \ldots, m, p \in L_k$$

**Proof.** Let $X_i = Y_{pi}^k + S_i + T_{pij}^k$. Obviously, $X_i$ is a normal random variable with mean $\mu_{Y_{pi}^k} + \mu_{S_i} + \mu_{T_{pij}^k}$ and variance $\sigma_{Y_{pi}^k}^2 + \sigma_{S_i}^2 + \sigma_{T_{pij}^k}^2$, $X_i$ is converted to the standard normal distribution by transformation function $X_i' = \dfrac{X_i - \mu_{Y_{pi}^k} - \mu_{S_i} - \mu_{T_{pij}^k}}{\sqrt{\sigma_{Y_{pi}^k}^2 + \sigma_{S_i}^2 + \sigma_{T_{pij}^k}^2}}$. Let $p_3 = P\left\{Y_{pi}^k + S_i + T_{pij}^k \leq h_j\right\}$. So, the probability $p_3$ is

$$p_3 = P\left\{Y_{pi}^k + S_i + T_{pij}^k \leq h_j\right\} = P\left\{X_i' \leq \frac{h_j - \mu_{Y_{pi}^k} - \mu_{S_i} - \mu_{T_{pij}^k}}{\sqrt{\sigma_{Y_{pi}^k}^2 + \sigma_{S_i}^2 + \sigma_{T_{pij}^k}^2}}\right\} = \Phi\left(\frac{h_j - \mu_{Y_{pi}^k} - \mu_{S_i} - \mu_{T_{pij}^k}}{\sqrt{\sigma_{Y_{pi}^k}^2 + \sigma_{S_i}^2 + \sigma_{T_{pij}^k}^2}}\right).$$

From $\mu_{Y_{pi}^k} = \mu_{A_{pi}^k} + \mu_{W_{pi}^k}$, $\sigma_{Y_{pi}^k}^2 = \sigma_{A_{pi}^k}^2 + \sigma_{W_{pi}^k}^2$, we have

$$P\left\{Y_{pi}^k + S_i + T_{pij}^k \leq h_j\right\} \geq 1 - \varepsilon_2 \Leftrightarrow \Phi\left(\frac{h_j - \mu_{A_{pi}^k} - \mu_{W_{pi}^k} - \mu_{S_i} - \mu_{T_{pij}^k}}{\sqrt{\sigma_{A_{pi}^k}^2 + \sigma_{W_{pi}^k}^2 + \sigma_{S_i}^2 + \sigma_{T_{pij}^k}^2}}\right) \geq 1 - \varepsilon_2,$$

$$i, j \in R, i \neq j, k = 1, \ldots, m, p \in L_k$$

□

Thus, the inequalities (7)–(9) in the model can be replaced by the inequalities (20), (23), and (26).

## 4. Proposed Population-Based Tabu Search Algorithm

The AMR scheduling problem is an NP-hard problem, which cannot be solved by any exact algorithm in polynomial time in general. The main solution method is heuristic algorithms. Given its stochastic nature and the fact that it combines several NP-hard subproblems, the SMMAMRRPTW is highly complex and difficult to solve exactly, even for relatively small instances. Tabu search (TS) algorithm is a meta-heuristic algorithm based on neighborhood search proposed by Glover [30] that is widely used to solve combinatorial optimization problems. It can solve a class of combinatorial optimization problems that involve stochastic factors. A genetic algorithm (GA) is a population-based meta-heuristic that can provide diverse solutions by searching the solution space parallelly [31].

In this section, we propose a population-based tabu search (PTS) algorithm that combines the tabu search, population mutation, and population crossover to solve SMMAMRRPTW. A repair operator is developed to repair the infeasible solutions. Algorithm 1 gives the overall process of the PTS algorithm. The population $Pop$, crossover probability $\lambda_1$, mutation probability $\lambda_2$, and the maximum number of iterations $N$ are the inputs for the algorithm. The solution $x_{best}$ is the output for the PTS algorithm. Line 1 initializes the operator weights $w_i$, tabu lists $B_i$, tabu tenure $t$, and iteration counter $ite$. Line 2 generates the initial chromosome $x_0$ of the problem through Algorithm 2. Line 3 takes the $x_0$ as the current optimal solution and the initial solution currently used for algorithm iteration. The main idea of the PTS algorithm is a while loop from lines 4 to 23. In each iteration, operator $i$ is selected through the roulette wheel procedure to generate neighborhood solution $N(x_{current})$, repair infeasible solutions in $N(x_{current})$ (line 5), and find the optimal neighborhood solution $x_{c\_best}$ from $N(x_{current})$ (line 6). Lines 7–11 compare the generated neighborhood solution $x_{c\_best}$ with the current optimal solution $x_{best}$ and update the current optimal solution, operator weights, and tabu table. The algorithm updates the weights $w_i$ of the operators every ten runs (line 12). Line 13 fills the population with percentage $r$ using $x_{best}$ and randomly fills the remaining population to achieve diversity. In lines 14–20, we first obtain the new population $Pop'$ using the crossover and mutation operators. Then, duplicate chromosomes are deleted in population $Pop'$, and finally, infeasible solutions are repaired in population $Pop'$. Population $Pop'$ is updated. When the algorithm reaches the maximum number of iterations $N$, the algorithm stops.

---

**Algorithm 1** Population-based tabu search (PTS) algorithm

**Input:** $Pop$, $\lambda_1$, $\lambda_2$, $N$.

**Output:** $x_{best}$

1. Initialize the operator weights $w_i = 1, i = 1, 2$, tabu lists $B_i, i = 1, 2$, tabu tenure $t$, and iteration counter $ite = 0$;
2. Generate initial chromosome $x_0$ through Greedy k-means algorithm; // **Algorithm 2**
3. $x_{current} \leftarrow x_0$, $x_{best} \leftarrow x_0$;
4. **While** $ite < N$ **do**
5. Use the roulette wheel procedure to select operators (2-opt and relocation) based on their weights $w_i$ to generate neighborhood solution $N(x_{current})$, and the depot insertion operator is used to repair infeasible solutions;
6. Find the optimal objective value $F(x_{c\_best})$, where $x_{c\_best} \in N(x_{current})$;
7. **If** $F(x_{c\_best}) < F(x_{best})$
8. $x_{best} \leftarrow x_{c\_best}$, $x_{current} \leftarrow x_{c\_best}$, $\rho_i \leftarrow \rho_i + 5$, and update the tabu list $B_i$;
9. **Else if**
10. $x_{current} \leftarrow x_{c\_best}$, $\rho_i \leftarrow \rho_i + 2$, and update the tabu list $B_i$;
11. **End if**
12. The weight $w_i = \frac{\rho_i}{\sum_{i=1,2} \rho_i}$ of each operator is updated after every ten iterations;
13. Fill the population with percentage $r$ using $x_{best}$ and randomly fill the remaining population;
14. **If** rand $< \lambda_1$ **then**
15. $Pop' \leftarrow crossover(Pop)$;
16. **End if**
17. **If** rand $< \lambda_2$ **then**
18. $Pop' \leftarrow mutation(Pop)$;
19. **End if**
20. Delete duplicate chromosomes, repair infeasible solutions in population;
21. The optimal chromosome in population $Pop'$ as the current solution $x_{current}$;
22. $ite \leftarrow ite + 1$;
23. **End while**
24. **Return** $x_{best}$

The structure of Section 4 is outlined as follows: in Section 4.1, for SMMAMRRPTW, we introduce the evaluation of the solution. Section 4.2 covers the creation of an initial solution. In Section 4.3, we introduce the framework of the tabu search algorithm. Sections 4.4 and 4.5 delve into the crossover operators and mutation operators. Lastly, Section 4.6 presents the termination condition for the PTS algorithm.

*4.1. Evaluation of Solution*

To facilitate the use of the PTS algorithm for solving SMMAMRRPTW, we convert multi-objective functions to a single objective. In terms of objective (2), to maximize the patient's satisfaction, the delayed service time of AMRs should be minimized, i.e.,

$$minmize\,(-f_2) = \sum_{k=1}^{m} \sum_{p \in L_k} \sum_{i \in R} E\left[\left(A_{pi}^k - h_i\right) \cdot z_{pi}^k\right]^+.$$

Hence, we can convert multi-objective functions to a single objective $minmize(f)$, where $f = f_1 + \xi_3(-f_2)$, $\xi_3$ is the penalty coefficient for delayed service time of AMRs.

*4.2. Initial Solution*

It is well known that initial solutions play an important role in the performance of the algorithm. High-quality initial solutions can make the algorithm converge quickly in a short computing time. K-means is widely used for its efficiency in cluster algorithms [32]. In k-means clustering, several data points are partitioned into clusters so that the within-cluster variances (squared Euclidean distances) are minimized.

In this paper, we employ a heuristic method of integrating Greedy insertion and k-means algorithm (Gk) for the PTS algorithm to construct the initial solution; Algorithm 2 provides the pseudocode. The maximum number of iterations $kiter$ of the Gk algorithm and the change in the threshold of cluster center point $\tau$ are the inputs for the algorithm. The solution $x_0$ is the output for the Gk algorithm. Randomly select $k$ requests as the initial centroids in Line 1. The Gk algorithm uses k-means to first group request points into distinct clusters; the main idea is a while loop from lines 2 to 6. In each iteration, the cluster, based on the geographical location of all patients, repeats multiple times with different initial centroids and selects the clustering result with the minimum loss function RSS, where RSS is the squared error, which is the Euclidean distance from each point to the centroid. And the requests are divided into $k$ clusters (line 7). Then, using the Greedy insertion algorithm, requests are assigned to the AMR for delivery within each cluster based on the time window size of the requests in each cluster. The subroutes $route_i$ in all clusters are connected to form an initial AMR route (lines 8–10). Finally, a repair operator is used to repair the initial route (line 11) to obtain the initial feasible solution $x_0$ of the PTS algorithm (line 12).

**Algorithm 2** Gk algorithm

**Input:** $kiter$, $\tau$
**Output:** solution $x_0$

1. Random select $k$ requests as the initial centroids;
2. **While** loss function RSS less than $1/\tau$ or not reach the number of iterations kiter **do**
3. Assign each request to its closest centroid to create clusters;
4. Compute loss function RSS;
5. Update centroids by taking means of request locations within each cluster;
6. **End while**
7. Divided into $k$ clusters $C = \{C_1, C_2, \cdots, C_k\}$;
8. **For** $C_i, i = 1, \cdots, k$ **do**
9. According to the lower limit of time window of request in each cluster $C_i$, connect requests to generate subroutes $route_i$;
10. **End for**
11. Using repair operators to obtain feasible solutions;
12. **Return** the feasible route as initial solution $x_0$

*4.3. Tabu Search*

Tabu search is the core composition of our algorithm. We select the 2-opt and relocation operator to the roulette wheel procedure to generate neighborhood solutions. The 2-opt operator, depicted in Figure 1, involves randomly selecting two nodes, 4 and 1 (as illustrated in Figure 1a), and employing the operator to interchange the segments between node 4 and node 1, resulting in a new solution, as depicted in Figure 1b. Similarly, the relocation operator, illustrated in Figure 2, begins by randomly selecting node 1 and node 4 (as shown in Figure 2a), followed by the insertion of node 1 in front of node 4. This process yields a new solution, as shown in Figure 2b.

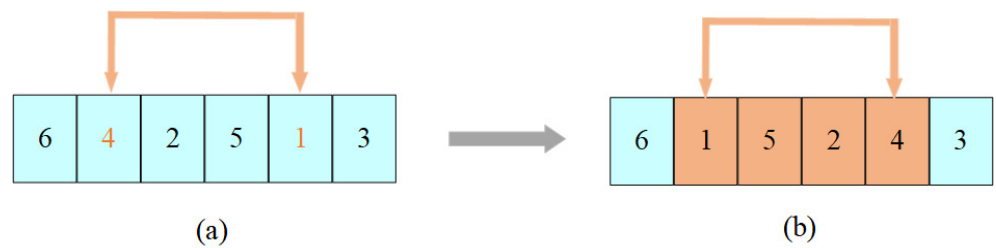

**Figure 1.** The 2-opt operator: (**a**) before using the 2-opt operator; (**b**) after using the 2-opt operator.

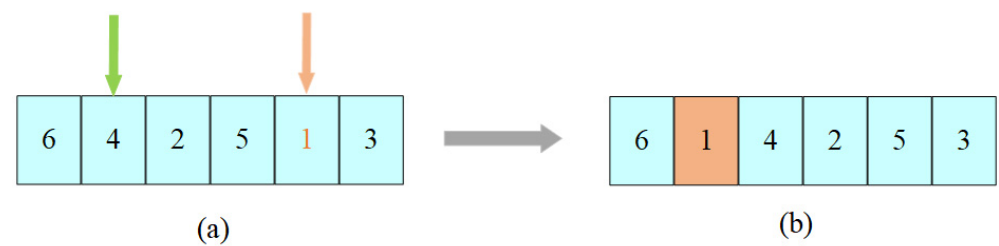

**Figure 2.** Relocation operator: (**a**) before using the relocation operator; (**b**) after using the relocation operator.

However, this operator does not always yield a feasible solution. To speed up the algorithm to find the optimal solution, we propose a repair operator (i.e., depot insertion operator) to repair the generated infeasible neighboring solution, as defined in Definition 1. The depot insertion operator is shown in Figure 3. For example, at node 5, the current route is infeasible (see Figure 3a). The depot insertion operator is used to insert a depot before the node where the route is infeasible; that is, insert node 0 before node 5 (see Figure 3b).

**Definition 1.** *The neighboring solution $x_{neighbor}$ is infeasible if the capacity constraints (6)–(8) or the time constraints (9)–(10) are not satisfied.*

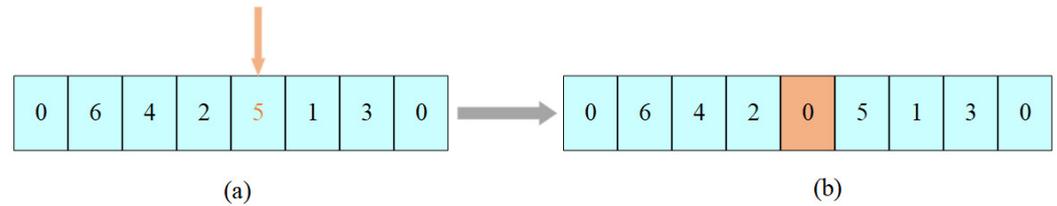

**Figure 3.** Depot insertion operator: (**a**) the current route is infeasible; (**b**) feasible route after repair.

For tabu lists, we use $B_i, \forall i = 1, 2$ to store the solutions produced by 2-opt and relocation to each pair of requests, where the tabu list $B_i$ is a $n \times n$ symmetric matrix, and $n$ is the total number of requests. For the operator $i$, the element $B_i(j_1, j_2)$ in the tabu list $B_i$ represents the tabu state of the neighborhood action $b_i(j_1, j_2)$ (i.e., tabu object), where $j_1, j_2 \in R$ are the two requests for 2-opt or relocation.

The tabu length in the PTS algorithm is used to record the number of iterations for removing tabu objects. As the PTS algorithm iterates, when the tabu length is zero, the tabu object $b_i(j_1, j_2)$ can be released from tabu table $B_i$.

In terms of the aspiration rule of the PTS algorithm, when a neighborhood solution is currently tabu but better than the current optimal solution, we revoke the tabu operation, and solution $x_{current}$ is selected to continue the algorithm iteration. Otherwise, the inferior solution shall be accepted. The suboptimal solution satisfying the tabu table is selected for iterations.

*4.4. Crossover Operators*

Crossover, often referred to as recombination, plays a vital role in evolutionary algorithms by combining the genetic information from two parents to produce offspring with hopefully better chromosomes. Essentially, crossover replaces partial structures of two-parent bodies to generate novel individuals. This operation not only retains the characteristics of high-quality individuals from the parent chromosomes but also maintains the diversity of solutions [33]. In this paper, before performing any operation on chromosome $P$, node 0 (this node representing depot) is removed and contained within it first. The crossover operation is shown in Figure 4, where fragments (4,2,5) are selected from the parent $P_1$, fragments (1,3,4) are selected from the parent $P_2$, and then fragments (4,2,5) and (1,3,4) are placed before the parent $P_2$ and $P_1$, respectively. The duplicate nodes are deleted, generating two new parents $P_1'$ and $P_2'$. Finally, the depot insertion operator is used to repair infeasibility in the offspring generated by the crossover operator.

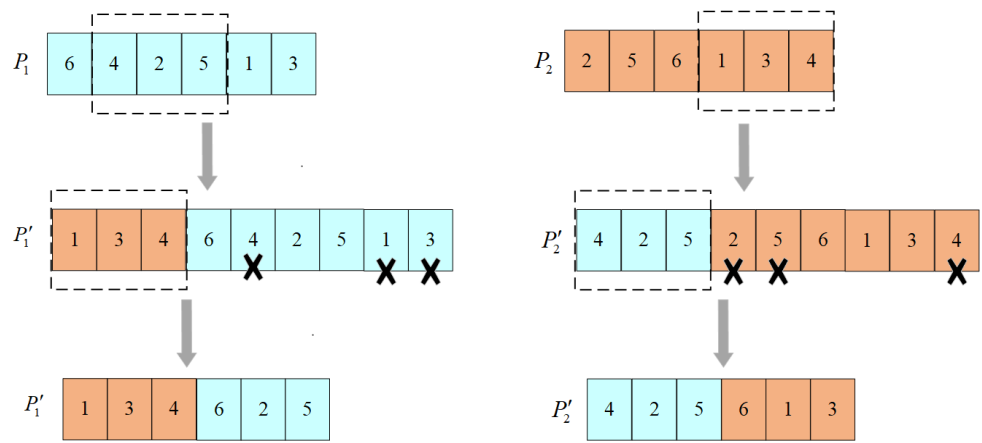

**Figure 4.** Crossover operator.

*4.5. Mutation Operators*

The mutation operator mainly imitates the gene mutation phenomenon in the process of biological evolution in nature and changes the individual gene values on the chromosome to obtain new chromosomes [34]. An excellent mutation operator can effectively avoid converging to the local optimum and generate high-quality solutions. The PTS algorithm proposed by this research selects the node swap mutation operator and the arc swap mutation operator to achieve diversity and widen the search space span.

In the arc swap mutation operator, two mutated gene segments are randomly generated, and the corresponding segments are exchanged. For example, in Figure 5, given a route of a chromosome (6,4,2,5,1,3) with randomly selected gene segments (6,4,2) and (1,3), the new route generated by this operator is (1,3,5,6,4,2). In the node swap mutation operator, two mutation gene bits are randomly generated, and the corresponding genes at the two positions are exchanged. In Figure 5, nodes 1 and 4 are selected and exchanged. The depot insertion operator is used to repair infeasible solutions of generated.

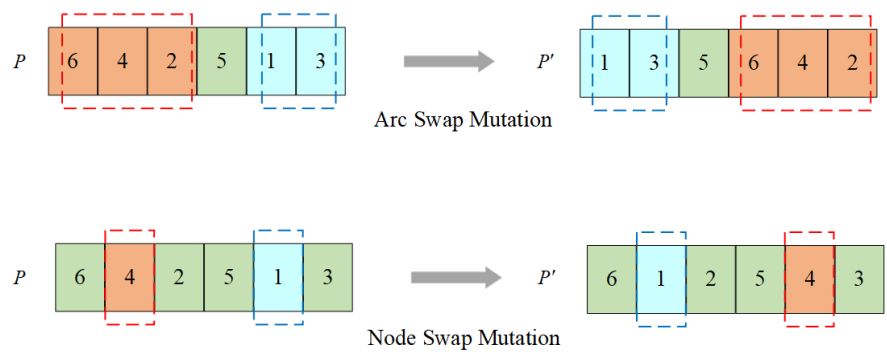

**Figure 5.** Swap operator.

*4.6. PTS Algorithm Stop Condition*

The algorithm stops after a predetermined maximum number of iterations $N$ and outputs the current optimal solution.

## 5. Computational Study

In this section, we first describe the benchmarks in the problem. Then, we report the parameter settings used in our experimental study. Moreover, the effectiveness of the developed algorithm is analyzed using comprehensive computational methods. Finally, we conducted relevant sensitivity analysis and management insights.

The algorithm is coded in MATLAB R2018b software. All experiments were performed on a PC with an Intel Core i5 CPU clocked at 2.40 GHz and 8 GB of memory running Windows 10 Professional 64-bit version.

*5.1. Dataset Generation*

The proposed model of SMMAMRRPTW involving stochastic demands and stochastic travel times has not been addressed in previous research. No comparative algorithm and benchmarks are available for comparison. Therefore, we use the well-known Solomon's dataset for VRPTW to generate our problem instances with appropriate modifications. The Solomon instances are grouped into three types. In type C, customers are clustering distributed. In type R, customers are randomly distributed, and in type RC, the customer distribution is mixed. Four instances are chosen from each type, and thus, there are a total of 12 instances tested, i.e., datasets C101, C102, C201, C202, R101, R102, R201, R202, RC101, RC102, RC201, and RC202.

We used the same vehicle capacity and customer and depot locations as in Solomon. We change the travel time of an edge $(i,j) \in A, i \neq j$ of an instance to a random variable $T_{pij}^k$, where the mean travel time $\mu_{T_{pij}^k}$ is the product of the Euclidean distance $d_{ij}$ between two points and the travel speed $v_r$, and the variance of $\sigma^2_{T_{pij}^k} = 0.2 \cdot \mu_{T_{pij}^k}$. The customer demands in Solomon as a mean ($\mu_{q_i}$, $i \in R$) for the stochastic demand in our problem, and the variance is $\sigma^2_{q_i} = 0.1 \cdot \mu_{q_i}$, $i \in R$. Furthermore, we assume that a service time per demand unit is 2 s ($a = 2s$) and the bias is 10 sec ($b = 10s$), which corresponds to approximately 12 s per demand unit. The confidence level for meeting the vehicle capacity $(1-\varepsilon_1)$ and time window $(1-\varepsilon_2)$ in SMMAMRRPTW is 0.95.

*5.2. Parameter Setting*

To determine the appropriate parameter settings, we used the Taguchi method to tune the parameters of the PTS algorithm. In the Taguchi method, orthogonal arrays are used to study a large number of decision variables with a small number of experiments. The Taguchi method offers orthogonal arrays to determine the best values of decision variables with a small number of experiments [35]. Taguchi method transforms the repetition data to another value, which is the signal-to-noise (S/N) ratio. Here, 'signal' denotes the response variable, and 'noise' denotes the standard deviation. Therefore, the S/N ratio indicates the amount of variation present in the response variable. The aim of the Taguchi method is to maximize the S/N ratio [35].

Table 2 shows the parameters of the PTS algorithm. For each parameter, three levels were considered. Ten different runs were performed on twelve instances, and the Minitab software 22 was used to provide the mean S/N ratio plot for the PTS algorithm, as shown in Figure 6. Based on the S/N ratio plot illustrated in Figure 6, the optimum level of all parameters was obtained (Table 2).

**Table 2.** Factor levels of the parameters of the PTS algorithm.

| Parameters | Factors | Parameter Levels | | | Optimum Level |
|---|---|---|---|---|---|
| | | Level 1 | Level 2 | Level 3 | |
| k-means iter | A | 20 | 30 | 40 | 30 |
| k-means threshold | B | 0.0001 | 0.0002 | 0.0003 | 0.0002 |
| Population size | C | 100 | 120 | 150 | 120 |
| Probability crossover | D | 0.7 | 0.8 | 0.9 | 0.7 |
| Probability mutation | E | 0.05 | 0.1 | 0.15 | 0.15 |
| Maximum number of iterations PTS | F | 50 | 200 | 500 | 50 |

Finally, the portion of the population to be filled was determined, and a set of runs with different percentages ($r$) for the six instances C101, C201, R101, R201, RC101, and RC201 was performed. Here, $r$ represents the population percentage filled by the current optimal solution of the TS algorithm in the PTS algorithm, while the remaining portion ($1-r$) is filled randomly to achieve diversity. Figure 7 shows the iteration diagrams of six instances at four different proportions $r$ (20%, 50%, 80%, and 100%). Figure 7 shows that using the current optimal solution generated by each iteration of the TS algorithm as a chromosome to fill the entire population ($r=100\%$) can make the PTS algorithm converge to a better solution.

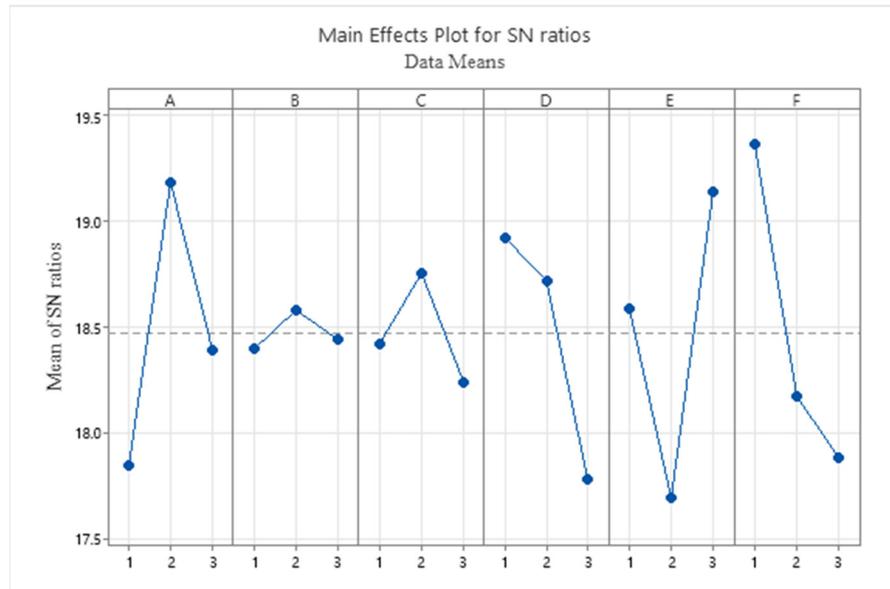

**Figure 6.** The mean S/N ratio plot for the PTS algorithm.

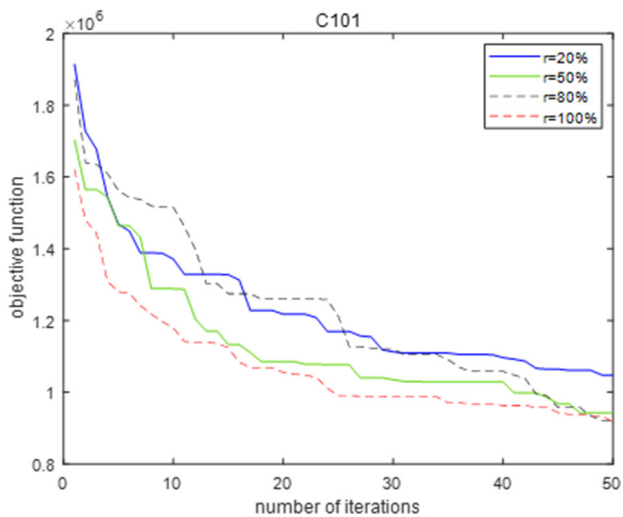
(a)

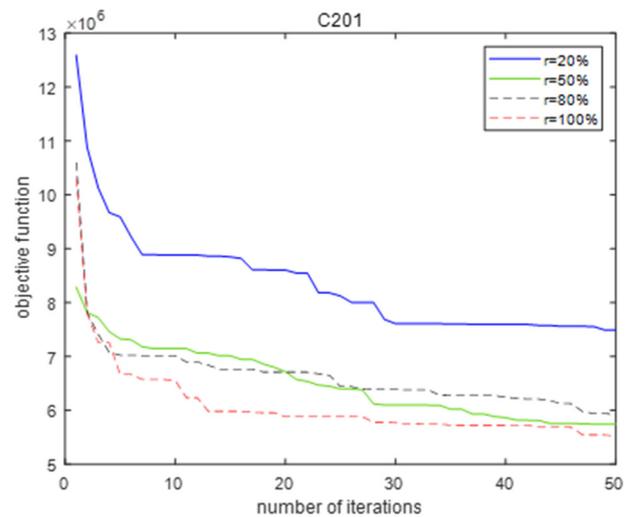
(b)

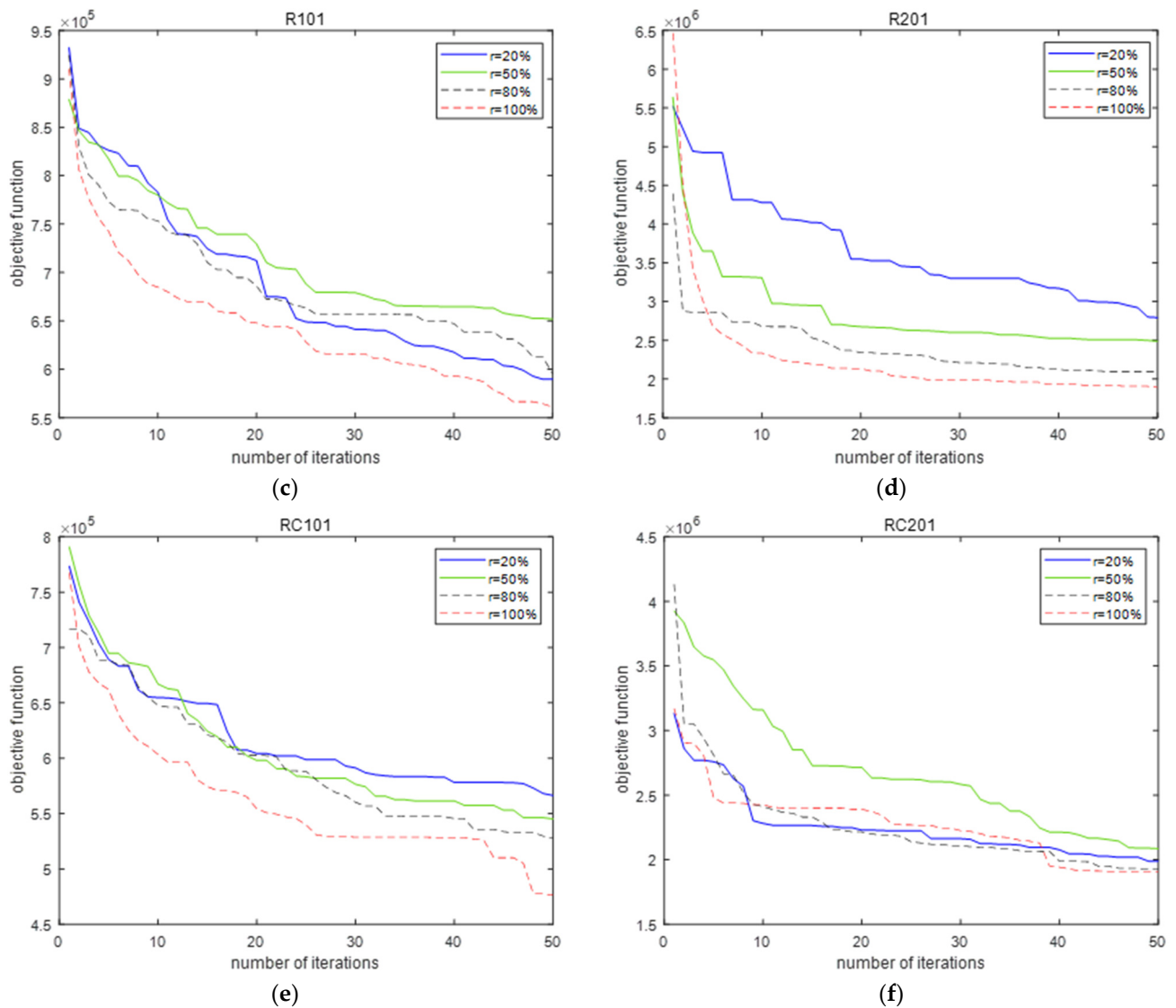

**Figure 7.** Experiment to determine the population setting of the PTS algorithm: (**a**) C101 instance; (**b**) C201 instance; (**c**) R101 instance; (**d**) R201 instance; (**e**) RC101 instance; (**f**) RC201 instance.

*5.3. Performance of the PTS Algorithm*

In this section, we conducted a series of experiments using 12 instances to assess the effectiveness and efficiency of the PTS algorithm and model. Specifically, the performance of the Gk algorithm for generating initial solutions was validated by comparing k-means and Greedy insertion algorithms. Furthermore, we compared the PTS with TS and GA algorithms and finally analyzed the service quality of patients.

We first compared Gk with the k-means and Greedy insertion algorithm on generating initial solutions. Each instance was run ten times, and the comparison results are described in Table 3, where $f_{ave}^{Gk}$ and $f_{ave}^{k}$ are the average objectives for ten runs, $f_{best}^{Gk}$ and $f_{best}^{k}$ are the optimal objectives for ten runs of the three algorithms, and $f^{G}$ is the objective of the Greedy insertion algorithm. It is worth noting that the Greedy insertion algorithm generates the same result every time it runs. The percentage improvement ($Imp_1$, $Imp_2$) of the average objective generated by the Gk algorithm relative to the average objectives generated by k-means and Greedy insertion algorithms, as defined by

$$Imp_1 = \frac{f_{ave}^k - f_{ave}^{Gk}}{f_{ave}^k} \times 100, \quad Imp_2 = \frac{f^G - f_{ave}^{Gk}}{f^G} \times 100.$$

We see that in most cases, Gk outperforms the k-means and Greedy insertion algorithm. Furthermore, we observe that the k-means method demonstrates favorable clustering outcomes for the type RC1 instances. This can be explained by the customer distribution characteristics and the width of time windows of RC1 instances. The k-means method can effectively cluster instances with mixed spatial distributions. The RC1 instance time window is relatively narrow, so the effectiveness of the Greedy insertion algorithm based on the width of time windows becomes less obvious. On the contrary, for the type C2 instances with clustering distribution and wide time windows, employing the Greedy insertion algorithm can obtain better initial solutions.

**Table 3.** Comparison results of Gk, k-means, and Greedy insertion algorithms for solving initial solutions.

| Instance | Gk | | k-Means | | | Greedy Insertion | |
|---|---|---|---|---|---|---|---|
| | $f_{ave}^{Gk}$ | $f_{best}^{Gk}$ | $f_{ave}^k$ | $f_{best}^k$ | $Imp_1$ | $f^G$ | $Imp_2$ |
| C101 | 2,045,910.08 | 1,863,415 | 2,610,782 | 2,417,649 | 21.63 | 2,617,170 | 21.82 |
| C102 | 1,254,631.38 | 1,107,483 | 1,655,803 | 1,434,495 | 24.22 | 1,870,543 | 32.92 |
| C201 | 12,085,514.7 | 8,324,882 | 14,809,985 | 11,851,861 | 18.39 | 7,814,317 | −54.65 |
| C202 | 9,612,167.96 | 8,365,007 | 12,044,963 | 9,850,724 | 20.19 | 9,311,346 | −3.23 |
| R101 | 938,552.54 | 884,402.9 | 959,881.3 | 830,621 | 2.22 | 1,166,671 | 19.55 |
| R102 | 742,706.32 | 686,893 | 797,058 | 739,009.9 | 6.81 | 953,384.5 | 22.09 |
| R201 | 5,093,775.34 | 3,129,344 | 7,378,546 | 6,652,661 | 30.96 | 5,823,723 | 12.53 |
| R202 | 4,614,830.48 | 2,115,929 | 6,124,538 | 5,632,099 | 24.65 | 6,014,803 | 23.27 |
| RC101 | 830,637.94 | 718,604 | 785,494.6 | 734,688.3 | −5.74 | 1,015,725 | 18.22 |
| RC102 | 693,004.58 | 628,846.4 | 666,737.6 | 576,412 | −3.93 | 929,352.6 | 25.43 |
| RC201 | 4,923,085.68 | 3,848,542 | 591,6827 | 5,416,624 | 16.79 | 5,317,144 | 7.41 |
| RC202 | 3,916,137.13 | 2,559,224 | 5,231,515 | 4,386,069 | 25.14 | 5,463,542 | 28.32 |
| Average | 3,895,912.85 | 2,852,714.40 | 4,915,177.43 | 4,210,242.97 | 20.73 | 4,024,810.2 | 3.20 |

We then compared PTS with the TS and GA algorithms on the SMMAMRRPTW problem. The comparative results are provided in Table 4, where we report the average objective ($f_{ave}^{TS}$, $f_{ave}^{GA}$, and $f_{ave}^{PTS}$), best objective ($f_{best}^{TS}$, $f_{best}^{GA}$, and $f_{best}^{PTS}$), Gap value ($G$), and the computation time (CPU) on ten runs of three algorithms. The results clearly show that our proposed PTS algorithm is markedly efficient. The relative percentage deviation ($G_1$, $G_2$) between the average objective value obtained by the PTS algorithm and the average objective value obtained by the TS and GA algorithm is calculated as follows, respectively.

$$G_1 = \frac{f_{ave}^{TS} - f_{ave}^{PTS}}{f_{ave}^{PTS}} \times 100, \quad G_2 = \frac{f_{ave}^{GA} - f_{ave}^{PTS}}{f_{ave}^{PTS}} \times 100.$$

A larger Gap value indicates better performance of the PTS algorithm. Indeed, the average objective and best objective in almost all the instances are found by the PTS algorithm. The average values of $f_{best}^{TS}$ and $f_{ave}^{TS}$ under the TS algorithm are 1,604,411 and 1,855,859, respectively, while the average values of $f_{best}^{PTS}$ and $f_{ave}^{PTS}$ under PTS are 1,481,096 and 1,707,478, respectively, demonstrating that PTS improves the average objective by $G_1 = 8.69\%$. This observation highlights the benefit of incorporating the population into the TS framework. Similarly, comparing PTS with the GA algorithm, the addition of the tabu list and aspiration rule improves the average objective by $G_2 = 35.6\%$. This result shows the benefits of incorporating the tabu list and aspiration rule into GA. In terms of runtime, we can see that the PTS algorithm proposed in this paper needs longer CPU computation time than the others, while the GA algorithm has the shortest computation time since each iteration of the PTS algorithm not only requires a tabu mechanism on the current neighborhood solution but also crossover and mutation of the population, which needs additional computation time.

**Table 4.** Comparison results of TS, GA, and PTS algorithms for solving SMMAMRRPTW.

| Instance | TS | | | GA | | | PTS | | | $G_1$ | $G_2$ |
|---|---|---|---|---|---|---|---|---|---|---|---|
| | $f_{best}^{TS}$ | $f_{ave}^{TS}$ | CPU (s) | $f_{best}^{GA}$ | $f_{ave}^{GA}$ | CPU (s) | $f_{best}^{PTS}$ | $f_{ave}^{PTS}$ | CPU (s) | | |
| C101 | 1,018,859 | 1,122,640 | 41.68 | 1,339,227 | 1,465,171 | 20.41 | 794,186.5 | 934,713.9 | 70.44 | 20.1 | 56.7 |
| C102 | 544,306.6 | 660,273.5 | 40.60 | 850,956.7 | 887,624.5 | 20.36 | 383,453.1 | 510,149.9 | 89.07 | 29.4 | 74.0 |
| C201 | 5,875,695 | 6,200,411 | 36.04 | 6,794,412 | 7,638,171 | 14.80 | 5,300,789 | 5,857,005 | 80.19 | 5.9 | 30.4 |
| C202 | 2,707,980 | 4,023,706 | 43.43 | 3,946,166 | 4,838,110 | 16.86 | 3,032,538 | 3,568,176 | 77.69 | 12.8 | 35.6 |
| R101 | 601,024.8 | 630,914.4 | 44.75 | 777,221.5 | 805,799.9 | 18.74 | 564,001.8 | 615,274.9 | 81.56 | 2.54 | 31.0 |
| R102 | 462,079 | 496,479.4 | 41.56 | 605,469.9 | 639,811.3 | 20.42 | 428,810.9 | 459,041.6 | 70.80 | 8.16 | 39.4 |
| R201 | 2,041,798 | 2,366,187 | 42.13 | 2,624,790 | 2,975,805 | 16.96 | 1,886,346 | 2,267,428 | 68.81 | 4.36 | 31.2 |
| R202 | 1,626,803 | 1,862,737 | 38.31 | 1,915,175 | 2,140,760 | 16.53 | 1,518,544 | 1,679,639 | 71.34 | 10.9 | 27.5 |
| RC101 | 535,131.9 | 572,251 | 48.61 | 652,803.2 | 709,288.6 | 20.72 | 484,415.7 | 548,133.4 | 72.90 | 4.4 | 29.4 |
| RC102 | 445,098.8 | 480,440.4 | 47.07 | 567,694.6 | 587,258.1 | 19.31 | 415,750.2 | 443,555.7 | 72.34 | 8.32 | 32.4 |
| RC201 | 1,888,989 | 2,147,014 | 43.4 | 2,682,160 | 2,991,024 | 17.10 | 1,689,098 | 2,106,756 | 72.36 | 1.91 | 42.0 |
| RC202 | 1,505,167 | 1,707,251 | 44.33 | 1,775,964 | 2,111,596 | 17.25 | 1,275,220 | 1,499,864 | 70.56 | 13.8 | 40.8 |
| Average | 1,604,411 | 1,855,859 | 42.66 | 2,044,337 | 2,315,868 | 18.29 | 1,481,096 | 1,707,478 | 74.84 | 8.69 | 35.6 |

All in all, the detailed analysis reveals that the PTS algorithm can find high-quality, objective value. Moreover, it is evident that incorporating population crossover and mutation operations into the TS algorithm significantly improves solution quality. This enhancement markedly boosts the scheduling performance of AMRs, particularly in complex routing scenarios.

*5.4. Sensitivity Analysis and Management Insights*

The CCP considers management decisions by controlling confidence levels when solving the SMMAMRRPTW problem. Our next experiment aims to analyze the impact of confidence level $1-\varepsilon_1$ of capacity constraint on the total delay service time (TDS) of AMRs and the number of AMRs. We compare the $1-\varepsilon_1 = 0.95$, $1-\varepsilon_1 = 0.9$, $1-\varepsilon_1 = 0.8$, $1-\varepsilon_1 = 0.7$, and $1-\varepsilon_1 = 0.6$, as shown in Figure 8.

From the observations in Figure 8a,b, we can see that increasing the confidence level tends to reduce the delayed service of AMRs, thereby improving the quality of service for patients. Especially for types C1, R1, and RC1 instances (Figure 8a), the impact of increasing the confidence level of capacity constraints on service quality is evident. Similarly, we can observe that increasing the confidence level leads to an increase in the use of AMR (see Figure 8c,d), which is more pronounced in types C1, R1, and RC1 instances (Figure 8c). This phenomenon can be explained by the smaller loading capacity of AMR in types C1, R1, and RC1 instances, while the larger capacity of AMR in types C2, R2, and RC2 instances. For instances with smaller capacity, increasing the confidence level on capacity constraints strengthens capacity constraints. Therefore, to satisfy capacity constraints, increasing the likelihood of AMR participating in services also increases accordingly, thereby reducing the cost of delayed services.

Based on the analysis above, we offer several management insights for hospitals regarding the scheduling and deployment of AMRs. Firstly, hospital managers should establish higher confidence levels to ensure prompt deliveries, particularly in critical areas such as medication and document transport. Secondly, adjusting confidence levels based on real demand and capacity constraints can significantly enhance operational efficiency.

This approach ensures optimal utilization of AMRs and improves overall performance. Lastly, for hospitals with frequent high service demands, investing in AMRs with larger capacities can be advantageous. This strategy reduces the frequency of trips and addresses capacity constraints with fewer AMRs, leading to long-term cost savings.

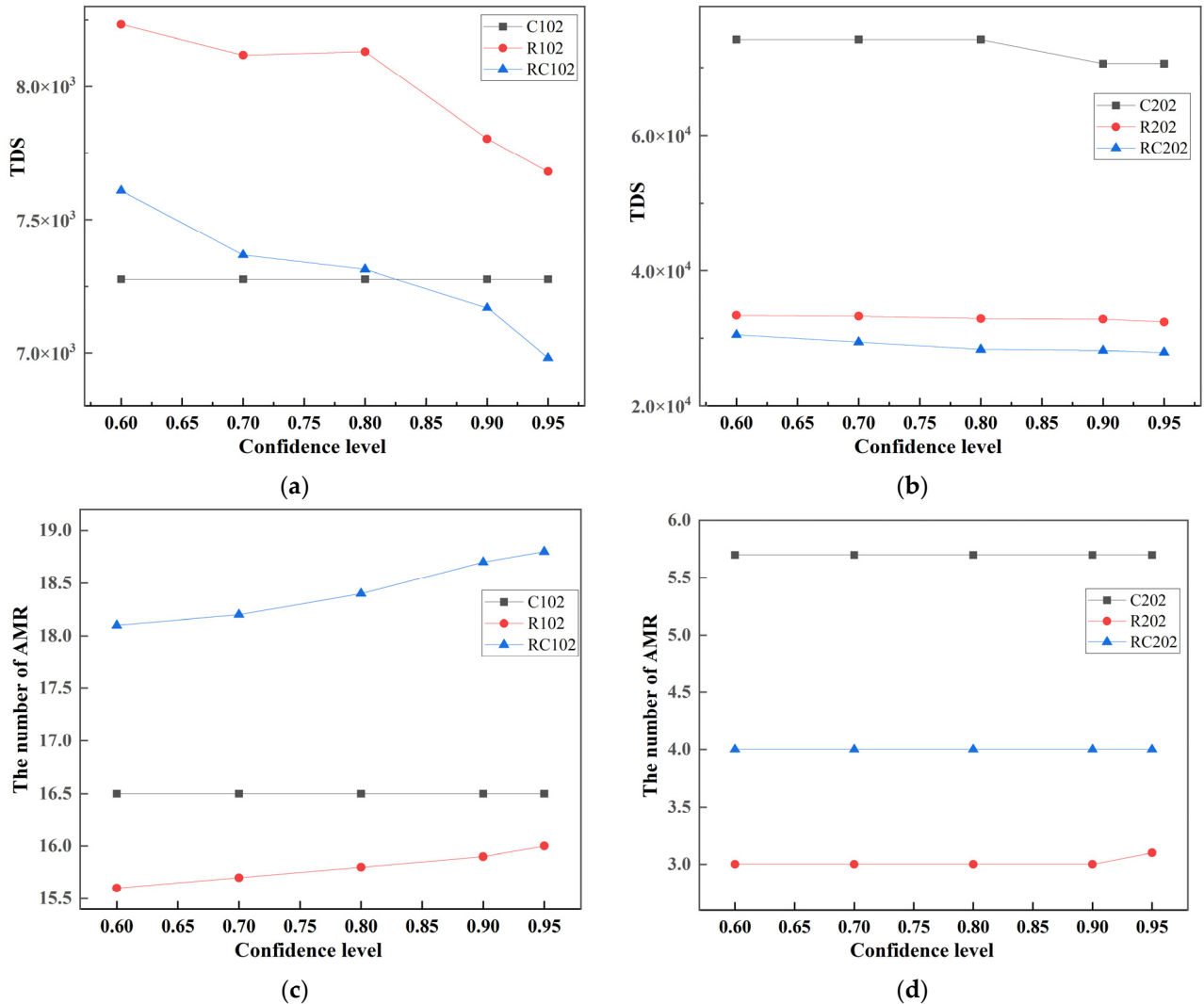

**Figure 8.** The impact of changing confidence levels: (**a**) the impact of confidence levels of types C1, R1, and RC1 instances on TDS; (**b**) the impact of confidence levels of types C2, R2, and RC2 instances on TDS; (**c**) the impact of confidence levels of types C1, R1, and RC1 instances on the number of AMR; (**d**) the impact of confidence levels of types C2, R2, and RC2 instances on the number of AMR.

## 6. Conclusions

In this paper, we studied the routing problem of autonomous mobile robots with multiple trips in the stochastic environment to minimize the total expected operating cost and maximize the total expected service quality of patients. The developed chance constraint models consider stochastic travel time, service time, and demands. By combining the tabu search (TS) algorithm with the genetic algorithm (GA), we designed the Population-based tabu search (PTS) algorithm, which is tailored to the characteristics of the stochastic multi-objective, multi-trip autonomous mobile robot routing problem with time windows (SMMAMRRPTW). By comparing the PTS algorithm with the TS and GA algorithms, we found that the PTS algorithm improves the objective by an average of $G_1 = 8.69\%$ compared to the TS algorithm, and by an average of $G_2 = 35.6\%$ compared to the GA algorithm. The proposed PTS algorithm incorporates several enhancement measures based on TS and GA, including initial solution generation, population generation, and neighborhood operators. These measures collectively and significantly improve the algorithm's performance. Our findings have significant practical implications. For hospital administrators, adjusting confidence levels based on actual demand and capacity constraints can significantly enhance operational efficiency, resulting in improved quality of service for patients and substantial operational cost savings. Furthermore, the methods and insights presented in this study can be extended to other logistics contexts beyond healthcare, such as warehousing and transportation.

Future research could expand our work in the following directions. Firstly, we focused on a single type of AMR for transportation services, but there is potential for multiple AMRs to transport various products. Thus, we are interested in extending our model to the collaborative transportation of multiple types of AMRs. Second, future research can explore the dynamic pick-up and delivery AMR scheduling problem with multi-depot. Lastly, future work may consider using reinforcement learning to analyze and solve these models.

**Author Contributions:** Conceptualization, L.C., N.Z. and K.W.; methodology, L.C., N.Z. and K.W.; software, L.C.; validation, L.C.; writing—original draft preparation, L.C.; writing—review and editing, L.C., N.Z. and K.W.; visualization, L.C.; supervision, L.C. and N.Z. All authors have read and agreed to the published version of the manuscript.

**Funding:** This research was supported in part by the State Key Laboratory of Industrial Control Technology (Grant No. ICT2021B51).

**Data Availability Statement:** The data used to support the findings of this paper are available from the corresponding author upon request.

**Conflicts of Interest:** The authors declare no conflicts of interest.